\documentclass[11pt,reqno]{article}
\usepackage{latexsym, amsmath, amssymb, a4, epsfig}
\usepackage{graphicx}

\newtheorem{theorem}{Theorem}[section]

\newcommand{\nm}{\noalign{\smallskip}}
\newcommand{\qed}{\hfill $\Box$}
\newcommand{\ds}{\displaystyle}

\newcommand{\p}{\partial}
\newcommand{\pd}[2]{\frac {\p #1}{\p #2}}
\newcommand{\eqnref}[1]{(\ref {#1})}

\newcommand{\Rbb}{\mathbb{R}}

\newcommand{\hatgrad}{\widehat{\nabla}}
\newcommand{\ckgrad}{\breve{\nabla}}
\renewcommand{\div}{\mbox{div}\,}

\newcommand{\la}{\langle}
\newcommand{\ra}{\rangle}

\def\Ba{{\bf a}}

\def\Bd{{\bf d}}
\def\Be{{\bf e}}
\def\Bf{{\bf f}}
\def\Bg{{\bf g}}

\def\Bt{{\bf t}}
\def\Bu{{\bf u}}
\def\Bv{{\bf v}}
\def\Bw{{\bf w}}
\def\Bx{{\bf x}}
\def\By{{\bf y}}

\def\BA{{\bf A}}
\def\BB{{\bf B}}

\def\BI{{\bf I}}

\newcommand{\Ga}{\alpha}
\newcommand{\Gb}{\beta}
\newcommand{\Gd}{\delta}
\newcommand{\Ge}{\epsilon}

\newcommand{\Gvf}{\varphi}

\newcommand{\Gs}{\sigma}

\newcommand{\GD}{\Delta}

\newcommand{\GG}{\Gamma}

\newcommand{\GO}{\Omega}

\newcommand{\beq}{\begin{equation}}
\newcommand{\eeq}{\end{equation}}

\newenvironment{dedication}
        {\vspace{1ex}\begin{quotation}\begin{center}\begin{em}}
        {\par\end{em}\end{center}\end{quotation}}

\numberwithin{equation}{section}
\numberwithin{figure}{section}

\begin{document}
\title{An over-determined boundary value problem arising from neutrally coated inclusions in three dimensions\thanks{\footnotesize This work is supported by the Korean Ministry of Education, Sciences and Technology through NRF grants Nos. 2010-0017532 and 2013R1A1A1A05009699, and by Grant-in-Aid
for Scientific Research (B) ($\sharp$ 26287020) of Japan Society for the Promotion of Science.}}

\author{Hyeonbae Kang\thanks{Department of Mathematics, Inha University, Incheon
402-751, Korea (hbkang, hdlee@inha.ac.kr).} \and Hyundae Lee\footnotemark[2] \and
Shigeru Sakaguchi\thanks{\footnotesize Research Center for Pure and Applied Mathematics,
Graduate School of  Information Sciences, Tohoku
University, Sendai, 980-8579,  Japan (sigersak@m.tohoku.ac.jp).}}

\maketitle

\begin{dedication}
Dedicated to the memory of Professor Kenjiro Okubo
\end{dedication}

\begin{abstract}
We consider the neutral inclusion problem in three dimensions which is to prove if a coated structure consisting of a core and a shell is neutral to all uniform fields, then the core and the shell must be concentric balls if the matrix is isotropic and confocal ellipsoids if the matrix is anisotropic. We first derive an over-determined boundary value problem in the shell of the neutral inclusion, and then prove in the isotropic case that if the over-determined problem admits a solution, then the core and the shell must be concentric balls. As a consequence it is proved that the structure is neutral to all uniform fields if and only if it consists of concentric balls provided that the coefficient of the core is larger than that of the shell.
\end{abstract}



%%%%%%%%%%%%%%%%%%%%%%%%%%%%%%%%%%%%%%%%%%%%%
\section{Introduction}
%%%%%%%%%%%%%%%%%%%%%%%%%%%%%%%%%%%%%%%%%%%%%

The purpose of this paper is to prove that the coated inclusions neutral to all uniform fields in the isotropic medium are concentric balls in three dimensions. The coated inclusion is depicted by $(D, \GO)$ where $D$ and $\GO$ are bounded  domains with Lipschitz boundaries in $\Rbb^d$ ($d=2,3$) such that $\overline{D} \subset \GO$. Here, $D$ represents the core and $\GO \setminus D$ the shell. The conductivity (or the dielectric constant) is $\Gs_c$ in the core and $\Gs_s$ in the shell ($\Gs_c \neq \Gs_s$). If the structure $(D, \GO)$ is inserted into the free space $\Rbb^d$ with conductivity $\Gs_m$ where there is a uniform field $-\nabla (\Ba \cdot \Bx)=-\Ba$ for some constant vector $\Ba$, then the field is perturbed in general. But for certain inclusions the field is not perturbed, in other words, the field does not recognize the existence of the inclusion. For example, the coated inclusion is made of concentric balls with specially chosen conductivities (confocal ellipsoids if $\Gs_m$ is anisotropic), one can see the uniform field is not perturbed. The inclusion with this property is called a neutral inclusion (or neutrally coated inclusion) and the neutral inclusion problem is to show that the inclusions of concentric balls (or confocal ellipsoids) are the only coated inclusions neutral to all uniform fields.

Let $\Gs$ denote the conductivity distribution of the medium so that
\beq
\Gs =
\begin{cases}
\Gs_c \quad&\mbox{in } D, \\
\Gs_s \quad&\mbox{in } \GO \setminus D, \\
\Gs_m \quad &\mbox{in } \Rbb^d \setminus \GO.
\end{cases}
\eeq
Here we assume that $\Gs_c$ and $\Gs_s$ are constants (or isotropic matrices), but $\Gs_m$ is allowed to be anisotropic symmetric matrix. We consider the following problem:
\beq\label{transgeneral}
\begin{cases}
\nabla \cdot \Gs \nabla u = 0 \quad \mbox{in } \Rbb^d, \\
u(\Bx)- \Ba \cdot \Bx = O(|\Bx|^{-2}) \quad \mbox{as } |\Bx| \to \infty,
\end{cases}
\eeq
where $\Ba$ is a constant vector. The term $u(\Bx)- \Ba \cdot \Bx$ depicts the perturbation of the potential due to insertion of the coated inclusion $(D,\GO)$. If the potential is not perturbed, namely,
\beq
u(\Bx)- \Ba \cdot \Bx \equiv 0 \quad\mbox{in } \Rbb^d \setminus \GO,
\eeq
the coated inclusion $(D, \GO)$ is said to be neutral to the field $\Ba$. If $(D, \GO)$ is neutral to the field $\Be_j$ for $j=1,\ldots, d$, where $\Be_j$ is the standard basis of $\Rbb^d$, then $(D, \GO)$ is neutral to all uniform fields.

Much interest in neutrally coated inclusions was aroused by the work of Hashin and Shtrikman \cite{HS} and Hashin \cite{hashine}. They showed that since insertion of neutral inclusions does not perturb the outside uniform field, the effective conductivity of the assemblage filled with coated inclusions of many different scales is $\Gs_m$. We refer to \cite{milton} for developments on neutral inclusions in relation to the theory of composites. Another interest in neutral inclusions has been aroused in relation to invisibility cloaking. The neutral inclusion is invisible from the probe by uniform fields as observed in \cite{kerker}.  Recently, the idea of neutrally coated inclusions has been extended to construct multi-coated circular structures which are neutral not only to uniform fields but also to fields of higher order up to $N$ for a given integer $N$ \cite{AKLL1}. It was proved there that the multi-coated structure combined with a transformation dramatically enhances the near cloaking of \cite{kohn1}. Cloaking by transformation optics was proposed in \cite{pendry} (and \cite{glu}).

As mentioned before, concentric balls (or disks) are neutral to all uniform fields by choosing $\Gs_c$, $\Gs_s$ and $\Gs_m$ properly ($\Gs_m$ is isotropic). Confocal ellipsoids (or ellipses) are also neutral to all uniform fields if $\Gs_m$ is anisotropic \cite{kerker} (see also section \ref{sec:ellipsoid}). Then a question arises naturally: are there any other shapes which are neutral to all uniform fields? In two dimensions there are no other shapes: if a coated inclusion $(D,\GO)$ is neutral to all uniform fields in two dimensions, then $D$ and $\GO$ are concentric disks (confocal ellipses if $\Gs_m$ is anisotropic). This is proved when $\Gs_c=0$ or $\infty$ in \cite{MS} and when $\Gs_c$ is finite in \cite{kl14}. In this paper we consider the neutral inclusion problem in three dimensions. We emphasize that the methods in \cite{kl14, MS} use powerful tools from complex analysis such as conformal mappings and harmonic conjugates, which cannot be applied to three dimensions. It is worth mentioning that there are many different shapes of coated inclusions neutral to a single uniform field as shown in two dimensions in \cite{JM, MS}.

We first show that if $(D,\GO)$ is neutral to all uniform fields in three dimensions and if $\Gs_c>\Gs_s$, then the following problem admits a solution:
\beq\label{free}
\left\{
\begin{array}{ll}
\ds \GD w= k \quad &\mbox{in } \GO \setminus \overline{D}, \\
\nm
\ds \nabla w = 0 \quad &\mbox{on } \p \GO , \\
\nm
\ds \nabla w(\Bx)= \BA\Bx +\Bd \quad &\mbox{on } \p D,
\end{array}
\right.
\eeq
where $k(\not=0)$ is a constant, $\BA$ is a symmetric matrix, and $\Bd$ is a constant vector. We emphasize that this is an over-determined problem because $\nabla w$ is prescribed on the boundaries. The problem, which is of independent interest, is to prove that if \eqnref{free} admits a solution in three dimensions, then $D$ and $\GO$ are confocal ellipsoids. If $D$ and $\GO$ are confocal ellipsoids, then \eqnref{free} admits a solution and $\BA$ should be either positive or negative-definite depending on the sign of $k$ (see section \ref{sec:ellipsoid}). So a part of the problem is to show that $\BA$ is either positive or negative-definite. In two dimensions it is proved in \cite{kl14} that if \eqnref{free} admits a solution then $D$ and $\GO$ are confocal ellipses (concentric disks if $\BA$ is isotropic). However, the proof there is based on the powerful result that there is a conformal mapping from $\GO \setminus \overline{D}$ onto an annulus. So it cannot be extended to three dimensions. The condition $\Gs_c>\Gs_s$, which is not natural, is required because of a technical reason for the derivation of \eqnref{free} in subsection \ref{sub22}. Even though we do not know how to do so, it is likely that the condition can be removed.

In this paper we solve the problem partially as the following theorem shows.

\begin{theorem}\label{conBall}
Let $D$ and $\GO$ be bounded  domains with Lipschitz boundaries in $\Rbb^3$ with $\overline{D} \subset \Omega$. Suppose that $\GO \setminus \overline{D}$ is connected.
If {\rm \eqnref{free}} admits a solution for $\BA=c\BI$ for some constant $c$ where $\BI$ is the identity matrix in three dimensions, then $D$ and $\GO$ are concentric balls.
\end{theorem}

As a consequence, we obtain the following theorem.

\begin{theorem}\label{thm:main}
Let $D$ and $\GO$ be bounded  domains with Lipschitz boundaries in $\Rbb^3$ with $\overline{D} \subset \Omega$.  Suppose that $\partial D$ is connected and $\mathbb R^3 \setminus \overline{D}$ is simply connected.  If $\Gs_m$ is isotropic, $\Gs_c>\Gs_s$ and $(D,\GO)$ is neutral to all uniform fields, then $D$ and $\GO$ are concentric balls.
\end{theorem}

This paper is organized as follows. In section \ref{sec:over} we show that if $(D,\GO)$ is neutral to all uniform fields then \eqnref{free} admits a solution. In section \ref{sec:ellipsoid} we construct a solution to \eqnref{free} when $D$ and $\GO$ are confocal ellipsoids. Section \ref{sec:proof} is to prove Theorem \ref{conBall}. In section \ref{sec:newton} we formulate the problem \eqnref{free} using Newtonian potentials and relate the problem with a known characterization of ellipsoids.

%%%%%%%%%%%%%%%%%%%%%%%%%%%%%%%%%%%%%%%%%%%%%
\section{Derivation of the over-determined problem}\label{sec:over}
%%%%%%%%%%%%%%%%%%%%%%%%%%%%%%%%%%%%%%%%%%%%%

In this section we derive \eqnref{free} out of the neutral inclusion problem. We will do so only in three dimensions since \eqnref{free} has been derived in two dimensions \cite{kl14}. We assume that $\partial D$ is connected and $\mathbb R^3 \setminus \overline{D}$ is simply connected.

Suppose, after diagonalization, that
\beq
\Gs_m = \mbox{diag} [\Gs_{m,1}, \Gs_{m,2}, \Gs_{m,3}].
\eeq
Let $u_j$, $j=1,2,3$, be the solution to
\beq\label{trans}
\left\{
\begin{array}{ll}
\nabla \cdot \Gs \nabla u_j = 0 \quad &\mbox{in } \Rbb^3,  \\
u_j(x)-x_j = O(|x|^{-2}) \quad &\mbox{as } |x| \to \infty.
\end{array}
\right.
\eeq
The structure being neutral to all three fields means that $u_j(x)-x_j=0$ in $\Rbb^3 \setminus \GO$ for $j=1,2,3$.
Let
\beq\label{wj}
w_j =  \frac{1}{\Gb_j} u_j
\eeq
where
$$
\Gb_j := \frac{\Gs_{m,j}}{\Gs_s}-1, \quad j=1,2,3,
$$
and $\Bw=(w_1,w_2,w_3)^T$ ($T$ for transpose). Set also
\beq\label{BB}
\BB= \mbox{diag} \, [1/\Gb_1, 1/\Gb_2, 1/\Gb_3].
\eeq
We will show the following:
\begin{itemize}
\item[(i)] $\nabla \Bw$ is symmetric and $\div \Bw$ is constant, and hence there is a function $\psi$ in $\overline{\GO}\setminus D$ such that
\beq\label{bwgrad}
\Bw=\nabla \psi\ \mbox{ and }\ \Delta \psi = \mbox{Tr}\, \BB +1 \quad\mbox{in }\GO \setminus \overline{D}.
\eeq
\item[(ii)] $\Bw(\Bx)=c_0\Bx+ \Bd$, $\Bx \in D$ for some constant $c_0$ and constant vector $\Bd$ (under the assumption that $\Gs_c > \Gs_s$).
\end{itemize}
We emphasize that it is in (ii) where the condition $\Gs_c > \Gs_s$ is required.

Once we have (i) and (ii), then we can show that \eqnref{free} has a solution. In fact, since $u_j=x_j$ on $\p\GO$, we have
$$
\nabla \psi(\Bx)= \BB \Bx \quad\mbox{on } \p\GO.
$$
Note that $\nabla \psi(\Bx)= c_0\Bx+ \Bd$ on $\p D$. Now define
\beq\label{definition_of_w}
w(\Bx):= \psi(\Bx)  - \frac{1}{2} \Bx \cdot \BB \Bx.
\eeq
Then $w$ satisfies \eqnref{free} with $k=1$ and $\BA=c_0 \BI - \BB$. We emphasize that if $\Gs_m$ is isotropic, so are $\BB$ and $\BA$.

%%%%%%%%%%%%%%%%%%%%%%%%%%%%%%%%%%%%%%%%%%%%%%%
\subsection{Proof of (i)}\label{shell}
%%%%%%%%%%%%%%%%%%%%%%%%%%%%%%%%%%%%%%%%%%%%%%%

Let us first deal with the case when $0< \Gs_c < \infty$. Denote by $\nu = (n_{1}, n_{2}, n_{3})^T$ the outward unit normal vector field to $\partial\Omega$ or $\partial D$. Note that the solution $u_j$ ($j=1,2,3$) to \eqnref{trans} satisfies the following transmission conditions on two interfaces :
\beq\label{trans1}
u_j|_+ - u_j|_- =0, \quad  \Gs_{m,j} \pd{u_j}{\nu}\Big|_+ - \Gs_s \pd{u_j}{\nu}\Big|_-=0 \quad \mbox{on } \p\GO
\eeq
and
\beq\label{trans2}
u_j|_+ - u_j|_- = 0, \quad \Gs_s \pd{u_j}{\nu}\Big|_+ - \Gs_c \pd{u_j}{\nu}\Big|_-=0 \quad \mbox{on } \p D
\eeq
where $+$ denotes the limit from outside and $-$ that from inside of $\GO$ or $D$. If $(D, \GO)$ is neutral to $x_j$, then $u_j(\Bx)-x_j=0$ in $\Rbb^3 \setminus \GO$, so we see from \eqnref{trans1} that
\beq\label{trans3}
u_j|_- = x_j, \quad \Gs_s \pd{u_j}{\nu}\Big|_-= \Gs_{m,j} n_j  \quad \mbox{on } \p\GO.
\eeq
In other words, $u_j$ is the solution to the following over-determined problem:
\beq\label{overdetermine}
\left\{
\begin{array}{ll}
\nabla \cdot \Gs \nabla u_j = 0 \quad &\mbox{in } \GO, \\
u_j = x_j, \ \  \ds \pd{u_j}{\nu} = \frac{\Gs_{m,j}}{\Gs_s} n_j  \quad &\mbox{on } \p\GO.
\end{array}
\right.
\eeq

Let $v_j \in C^2(\overline{\GO})$. Then we see from the divergence theorem and \eqnref{trans2} that
\begin{align*}
\int_{\p\GO} \pd{u_j}{\nu}\Big|_- v_j - u_j \pd{v_j}{\nu}
&= - \int_{\GO \setminus D} u_j \Delta v_j + \int_{\p D} \pd{u_j}{\nu}\Big|_+ v_j - u_j \pd{v_j}{\nu} \\
&= - \int_{\GO \setminus D} u_j \Delta v_j + \left( \frac{\Gs_c}{\Gs_s} -1 \right) \int_{\p D} \pd{u_j}{\nu}\Big|_- v_j + \int_{\p D} \pd{u_j}{\nu}\Big|_- v_j - u_j \pd{v_j}{\nu} \\
&= - \int_{\GO \setminus D} u_j \Delta v_j + \left( \frac{\Gs_c}{\Gs_s} -1 \right) \int_{D} \nabla u_j \cdot \nabla v_j - \int_{D} u_j \Delta v_j \\
&= - \int_{\GO} u_j \Delta v_j + \left( \frac{\Gs_c}{\Gs_s} -1 \right) \int_{D} \nabla u_j \cdot \nabla v_j .
\end{align*}
On the other hand, we see from \eqnref{trans3} that
\begin{align*}
\int_{\p\GO} \pd{u_j}{\nu}\Big|_- v_j - u_j \pd{v_j}{\nu}
&= \int_{\p\GO} \frac{\Gs_{m,j}}{\Gs_s} n_j v_j - y_j \pd{v_j}{\nu} \\
&= \left( \frac{\Gs_{m,j}}{\Gs_s} -1 \right) \int_{\GO} \pd{v_j}{y_j} - \int_{\GO} y_j \Delta v_j.
\end{align*}
Equating two identities above we obtain
\beq\label{weakform}
\int_{\GO} (y_j-u_j) \Delta v_j + \Ga \int_{D} \nabla u_j \cdot \nabla v_j - \Gb_j \int_{\GO} \pd{v_j}{y_j}=0 , \quad j=1,2,3
\eeq
for $v_j \in C^2(\overline{\GO})$, where $\Ga$ and $\Gb_j$ are defined for ease of notation to be
\beq
\Ga= \frac{\Gs_c}{\Gs_s} -1 \quad\mbox{and} \quad \Gb_j= \frac{\Gs_{m,j}}{\Gs_s}-1.
\eeq
Let $w_j$ be defined by $w_j := \frac{1}{\Gb_j} u_j$ as in \eqnref{wj}.
Then \eqnref{weakform} can be rephrased as
\beq\label{weakform11}
\int_{\GO} (\frac{1}{\Gb_j} y_j-w_j) \GD v_j + \Ga \int_{D} \nabla w_j \cdot \nabla v_j - \int_{\GO} \pd{v_j}{y_j}=0 , \quad j=1,2,3.
\eeq
Summing \eqnref{weakform11} over $j=1,2,3$ we have
$$
\int_{\GO} \sum_{j=1}^3 (\frac{1}{\Gb_j}y_j-w_j) \GD v_j + \Ga \int_{D} \sum_{j=1}^3 \nabla w_j \cdot \nabla v_j - \int_{\GO} \sum_{j=1}^3  \pd{v_j}{y_j}=0
$$
for $v_j \in C^2(\overline{\GO})$. If we use vector notation $\Bw=(w_1,w_2, w_3)^T$ and $\Bv=(v_1,v_2,v_3)^T$ ($T$ for transpose), then the above identity can be rewritten as
\beq\label{weakform2}
\int_{\GO} (\BB\By-\Bw) \cdot \Delta \Bv + \Ga \int_{D} \nabla \Bw : \nabla \Bv - \int_{\GO} \div \Bv =0.
\eeq
Here and afterwards $\BA:\BB$ denote the contraction of two matrices $\BA$ and $\BB$,
{\it i.e.}, $\BA:\BB=\sum a_{ij}b_{ij}=\textrm{Tr}(\BA^{T}\BB)$.

Let $\GG$ be the fundamental solution of the Laplace operator in $\Rbb^3$, {\it i.e.},
\beq
\GG(\Bx) := - \frac{1}{4 \pi |\Bx|}, \quad \Bx \neq 0.
\eeq
Let $v_j(\By) =\GG(\Bx-\By)$ for a fixed $\Bx \in \GO$.
Since $\Delta v_j(\By) = \Gd(\Bx-\By)$, by applying the divergence theorem over $\GO \setminus B_\Ge(\Bx)$ for sufficiently small $\Ge$ ($B_\Ge(\Bx)$ is the ball of radius $\Ge$ centered at $\Bx$) we see from \eqnref{weakform11} that
\beq\label{repre}
w_j (\Bx) = \frac{1}{\Gb_j} x_j + \Ga \int_{D} \nabla w_j(\By) \cdot \nabla_\By \GG(\Bx-\By)d\By + \pd{}{x_j} N_\GO(\Bx) , \quad \Bx \in \GO, \quad j=1,2,3,
\eeq
where $N_\GO$ is the Newtonian potential on a domain $\GO$, {\it i.e.},
\beq
N_\GO(\Bx) := \int_\GO \GG(\Bx-\By) d\By, \quad \Bx \in \Rbb^3.
\eeq
Let
$$
f_j(\Bx) := \int_{D} \nabla w_j(\By) \cdot \nabla_\By \GG(\Bx-\By) d\By, \quad j=1,2,3,
$$
and let $\Bf= (f_1, f_2, f_3)^T$. Note that $f_j$ is harmonic in $\Rbb^3 \setminus \overline{D}$, and \eqnref{repre} can be rewritten as
\beq\label{ujy}
\Bw(\Bx) = \Ga \Bf(\Bx) + \nabla \left( \frac{1}{2} \Bx\cdot \BB\Bx + N_\GO(\Bx) \right) , \quad \Bx \in \GO, \quad j=1,2,3.
\eeq

For any fixed $\Bx \in \Rbb^3 \setminus \overline{\GO}$, let
$$
v_j(\By)= \pd{}{x_j} \GG(\Bx-\By), \quad j=1,2,3.
$$
Then $\div \Bv (\By)= -\GD_\By \GG(\Bx-\By)=0$ and $\GD \Bv(\By)=0$ for $\By \in \GO$. So we see from \eqnref{weakform2} that
$$
\int_D \nabla \Bw : \nabla \Bv =0,
$$
and hence
\beq\label{sumfj}
\div \Bf(\Bx) = \int_{D} \sum_j \nabla w_j(\By) \cdot \nabla \pd{}{x_j} \GG(\Bx-\By) d\By= \int_D \nabla \Bw : \nabla \Bv =0, \quad \Bx \in \Rbb^3 \setminus \overline{\GO}.
\eeq
Since $f_j$ is harmonic in $\Rbb^3 \setminus \overline{D}$, \eqnref{sumfj} holds for all $\Bx \in \Rbb^3 \setminus \overline{D}$.

Again fix $\Bx \in \Rbb^3 \setminus \overline{\GO}$. Let $\{i,j,k \}$ be a permutation of $\{ 1,2,3 \}$ and let
$$
v_i(\By)= \pd{}{x_j} \GG(\Bx-\By), \quad v_j(\By)= -\pd{}{x_i} \GG(\Bx-\By), \quad v_k=0, \quad \By \in \GO.
$$
Then, $\GD \Bv=0$ and $\div \Bv=0$ in $\GO$. So we have from \eqnref{weakform2}
$$
\int_{D} \nabla w_i(\By) \cdot \nabla \pd{}{x_j} \GG(\Bx-\By) d\By - \int_{D} \nabla w_j(\By) \cdot \nabla \pd{}{x_i} \GG(\Bx-\By) d\By=0,
$$
which implies that
\beq
\p_i f_j(\Bx) = \p_j f_i(\Bx)
\eeq
for all $\Bx \in \Rbb^3 \setminus \overline{\GO}$ and hence for all $\Bx \in \Rbb^3 \setminus \overline{D}$. Moreover,  since $\Rbb^3 \setminus \overline{D}$ is simply connected, by the Stokes theorem
%since $f_j(\Bx) = O(|\Bx|^{-2})$ as $|\Bx| \to \infty$,
there is $\Gvf$ such that
\beq\label{fjy}
\Bf(\Bx) = \nabla \Gvf(\Bx), \quad \Bx \in \Rbb^3 \setminus \overline{D}.
\eeq
Because of \eqnref{sumfj}, we have
\beq\label{Deltag}
\GD \Gvf(\Bx)=0 , \quad \Bx \in \Rbb^3 \setminus \overline{D}.
\eeq

Let
\beq
\psi(\Bx) = \Ga \Gvf(\Bx) + \frac{1}{2} \Bx \cdot \BB\Bx +  N_\GO(\Bx)  , \quad \Bx \in \GO \setminus \overline{D}.
\eeq
Then, we have from \eqnref{ujy} and \eqnref{fjy}
\beq\label{ujy2}
\Bw (\Bx) = \nabla \psi(\Bx) , \quad \Bx \in \GO\setminus \overline{D}, \quad j=1,2,3.
\eeq
Since $\GD N_\GO(\Bx)=1$ for $\Bx \in \GO$, we have from \eqnref{Deltag} that
\beq
\GD \psi(\Bx)= \mbox{Tr}\, \BB + 1, \quad \Bx \in \GO\setminus \overline{D}.
\eeq
 So far we have shown that $\nabla \Bw$ is symmetric, $\div \Bw$ is constant, and \eqnref{bwgrad} holds  when $\Gs_c$ is finite.

We now assume that $\Gs_c=0$. In this case the problem \eqnref{overdetermine} becomes
\beq\label{overzero}
\left\{
\begin{array}{ll}
\GD u_j = 0 \quad &\mbox{in } \GO \setminus \overline{D}, \\
\nm
\ds \pd{u_j}{\nu} =0  \quad &\mbox{on } \p D, \\
\nm
u_j = x_j, \ \  \ds \pd{u_j}{\nu} = \frac{\Gs_{m,j}}{\Gs_s} n_j  \quad &\mbox{on } \p\GO.
\end{array}
\right.
\eeq
So, we see in a way similar to \eqnref{weakform} that
\beq\label{weakzero}
\int_{\GO} y_j \GD v_j - \int_{\GO \setminus D} u_j \GD v_j - \int_{\p D} u_j \pd{v_j}{\nu} - \Gb_j \int_{\GO} \pd{v_j}{y_j}=0
\eeq
for all $v_j \in C^2(\overline{\GO})$. So we obtain a representation of the solution similar to \eqnref{repre}:
\beq\label{reprezero}
w_j (\Bx) = \frac{1}{\Gb_j} x_j - \int_{\p D} w_j(\By) \pd{}{\nu}\GG(\Bx-\By)d\sigma(\By) + \pd{}{x_j} N_\GO(\Bx) , \quad \Bx \in \GO \setminus D.
\eeq
So, we infer in the exactly same way as in the previous sections that $\nabla\Bw$ is symmetric and $\div \Bw$ is constant, and there is a function $\psi$ such that \eqnref{bwgrad} holds.

Suppose that $\Gs_c=\infty$. In this case the problem \eqnref{overdetermine} becomes
\beq\label{overinfty}
\left\{
\begin{array}{ll}
\GD u_j = 0 \quad &\mbox{in } \GO \setminus \overline{D}, \\
u_j = \gamma_j \ (\mbox{constant}) \quad &\mbox{on } \p D, \\
\nm
u_j = x_j, \ \  \ds \pd{u_j}{\nu} = \frac{\Gs_{m,j}}{\Gs_s} n_j  \quad &\mbox{on } \p\GO.
\end{array}
\right.
\eeq
The constant $\gamma_j$ is determined by the condition
$$
\int_{\p D} \pd{u_j}{\nu} \Big|_+ =0.
$$
We then obtain similarly to \eqnref{weakform}
\beq\label{weakinfty}
\int_{\GO} y_j \GD v_j - \int_{\GO \setminus D} u_j \GD v_j + \int_{\p D} \pd{u_j}{\nu} v_j - \gamma_j \int_{D} \GD v_j - \Gb_j \int_{\GO} \pd{v_j}{y_j}=0
\eeq
for all $v_j \in C^2(\overline{\GO})$. We then obtain a representation of the solution similar to \eqnref{repre}:
\beq\label{repreinfty}
w_j (\Bx) = \frac{1}{\Gb_j} x_j + \int_{\p D} \pd{w_j}{\nu}(\By) \GG(\Bx-\By)d\sigma(\By) + \pd{}{x_j} N_\GO(\Bx) , \quad \Bx \in \GO \setminus D.
\eeq
So, we infer that $\nabla \Bw$ is symmetric, $\div \Bw$ is constant, and there is a function $\psi$ such that \eqnref{bwgrad} holds.

%%%%%%%%%%%%%%%%%%%%%%%%%%%%%%%%%%
\subsection{Proof of (ii)}\label{sub22}
%%%%%%%%%%%%%%%%%%%%%%%%%%%%%%%%%%

The transmission conditions \eqnref{trans2} on $\p D$ can be rephrased as
\beq\label{transvec}
\Bw|_+= \Bw|_-,  \quad \Gs_s \nabla \Bw|_+ \nu = \Gs_c \nabla \Bw|_- \nu.
\eeq
Let $\Bt_1$ and $\Bt_2$ be two orthonormal tangent vector fields to $\p D$.
Then, we have
$$
(\div \Bw)_- = \la (\nabla \Bw)_- \nu, \nu \ra + \la (\nabla \Bw)_- \Bt_1, \Bt_1 \ra +  \la (\nabla \Bw)_- \Bt_2, \Bt_2 \ra,
$$
and
$$
(\div \Bw)_+ = \la (\nabla \Bw)_+ \nu, \nu \ra + \la (\nabla \Bw)_+ \Bt_1, \Bt_1 \ra +  \la (\nabla \Bw)_+ \Bt_2, \Bt_2 \ra.
$$
Here $(\div \Bw)_-$ denotes the limit of $\div \Bw$ to $\p D$ from inside $D$, and $(\div \Bw)_+$ denotes that from outside $D$. Since
$$
\la (\nabla \Bw)_- \Bt_j, \Bt_j \ra = \la (\nabla \Bw)_+ \Bt_j, \Bt_j \ra, \quad j=1,2,
$$
we have
$$
(\div \Bw)_- - (\div \Bw)_+ = \la (\nabla \Bw)_- \nu, \nu \ra - \la (\nabla \Bw)_+ \nu, \nu \ra.
$$
It then follows from the second identity in \eqnref{transvec} that
\beq\label{transvec2}
\left\la \big( (\nabla \Bw)_-^T - \frac{\Gs_c}{\Gs_s} (\nabla \Bw)_- \big) \nu, \nu \right\ra = (\div \Bw)_- - (\div \Bw)_+ .
\eeq
On the other hand, since $(\nabla \Bw)_+$ is symmetric, we obtain
\begin{align}
\left\la \big( (\nabla \Bw)_-^T - \frac{\Gs_c}{\Gs_s} (\nabla \Bw)_- \big) \nu, \Bt_j \right\ra
& = \la \nu, (\nabla \Bw)_- \Bt_j \ra - \frac{\Gs_c}{\Gs_s} \la  (\nabla \Bw)_- \nu, \Bt_j \ra \nonumber \\
& = \la \nu, (\nabla \Bw)_+ \Bt_j \ra - \la  (\nabla \Bw)_+ \nu, \Bt_j \ra =0. \label{transvec3}
\end{align}
We then infer from \eqnref{transvec2} and \eqnref{transvec3} that
\beq\label{transvec4}
\left( (\nabla \Bw)_-^T - \frac{\Gs_c}{\Gs_s} (\nabla \Bw)_- \right) \nu = (\div \Bw)_- \nu - (\div \Bw)_+ \nu.
\eeq

Recall that $\div \Bw$ is constant in $\GO \setminus D$. Let
\beq
\Bv(\Bx)= \Bw(\Bx) - \frac{(\div \Bw)_+}{2+ \frac{\Gs_c}{\Gs_s}} \Bx, \quad \Bx \in D.
\eeq
Then one can see from \eqnref{transvec4} that
\beq\label{zerobdry}
\left( (\nabla \Bv)^T - \frac{\Gs_c}{\Gs_s} (\nabla \Bv) \right) \nu - (\div \Bv) \nu =0 \quad \mbox{on } \p D.
\eeq
Let $\Bg$ be a smooth vector field on $\overline{D}$.
It follows from \eqnref{zerobdry} and the divergence theorem that
\begin{align*}
0 &= \int_{\p D} \nu \cdot \left( (\nabla \Bv)\Bg - \frac{\Gs_c}{\Gs_s} (\nabla \Bv)^T \Bg - (\div \Bv) \Bg \right) d\Gs \\
&= \int_{D} \div \left( (\nabla \Bv)\Bg - \frac{\Gs_c}{\Gs_s} (\nabla \Bv)^T \Bg - (\div \Bv) \Bg \right) d\Bx.
\end{align*}
One can easily show that
$$
\div \left( (\nabla \Bv)\Bg - \frac{\Gs_c}{\Gs_s} (\nabla \Bv)^T \Bg - (\div \Bv) \Bg \right) =
\nabla \Bv^T : \nabla \Bg - \frac{\Gs_c}{\Gs_s} \nabla \Bv: \nabla \Bg - (\div \Bv)(\div \Bg),
$$
and so we obtain
\beq
\int_D \nabla \Bv^T : \nabla \Bg - \frac{\Gs_c}{\Gs_s} \nabla \Bv: \nabla \Bg - (\div \Bv)(\div \Bg)=0.
\eeq
Using notation
$$
\hatgrad \Bv :=\frac{1}{2}(\nabla \Bv +\nabla \Bv^{T}) \quad\mbox{and}\quad \ckgrad \Bv :=\frac{1}{2}(\nabla \Bv -\nabla \Bv^{T}),
$$
it can be rewritten as
\beq
(1- \frac{\Gs_c}{\Gs_s}) \int_D \hatgrad \Bv: \hatgrad\Bg - (1+\frac{\Gs_c}{\Gs_s}) \int_D \ckgrad \Bv: \ckgrad\Bg - \int_D (\div \Bv)(\div \Bg)=0.
\eeq

If $\Gs_c > \Gs_s$, then we take $\Bg=\Bv$ so that
\beq\label{239}
(1- \frac{\Gs_c}{\Gs_s}) \int_D |\hatgrad \Bv|^2 - (1+\frac{\Gs_c}{\Gs_s}) \int_D |\ckgrad \Bv|^2 - \int_D (\div \Bv)^2=0.
\eeq
Thus, we infer that $\Bv$ is constant in $D$ and hence
\beq
\Bw(\Bx) = \frac{(\div \Bw)_+}{2+ \frac{\Gs_c}{\Gs_s}} \Bx + \mbox{a constant vector}, \quad \Bx \in D.
\eeq
If $\Gs_c = \infty$, then $\Bu$ is constant on $\p D$, and hence $(\nabla \Bw)\Bt=0$ on $\p D$ for any tangential vector $\Bt$ to $\p D$. Since $\nabla \Bw$ is symmetric and $\div \Bw$ is constant, it implies that
$$
(\nabla \Bw)\nu = c_0\nu  \quad \mbox{on } \p D
$$
for some constant $c_0$. So, we can see that (ii) holds.

%%%%%%%%%%%%%%%%%%%%%%%%%%%%%%%%%%%%%%%%%%%%%%%%%%%%%%%%%
\section{Existence of solutions on confocal ellipsoids}\label{sec:ellipsoid}
%%%%%%%%%%%%%%%%%%%%%%%%%%%%%%%%%%%%%%%%%%%%%%%%%%%%%%%%%

We first mention that the solution $w$ to \eqnref{free}  is unique in the sense that if $w_1$ and $w_2$ are two solutions (with different $k,\, \BA$'s, and $\Bd$'s), then $w_1=C w_2 +E$ for some constants $C$ and $E$.  In fact, if $w_j$ is a solution to \eqnref{free} with $k=k_j \neq 0$, $\BA=\BA_j$ and $\Bd=\Bd_{j}$ ($j=1,2$), then $w=w_1-\frac{k_1}{k_2} w_2$ satisfies $\GD w= 0$ in $\GO \setminus \overline{D}$ and $\nabla w = 0$ on $\p\GO$, so we have that $w$ must be a constant.

We now construct a solution to \eqnref{free} when $D$ and $\GO$ are confocal ellipsoids. To do so, assume that $\p D$ is given by
\beq \label{defb3}
\frac{x^2_1}{c^2_1}+\frac{x^2_2}{c^2_2}+\frac{x^2_3}{c^2_3}=1.
\eeq
We then use the confocal ellipsoidal coordinates
$\rho, \ \mu,\ \xi$ such that
\begin{align*}
\frac{x^2_1}{c^2_1+\rho}+\frac{x^2_2}{c^2_2+\rho}+\frac{x^2_3}{c^2_3+\rho}=1,\\
\nm
\frac{x^2_1}{c^2_1+\mu}+\frac{x^2_2}{c^2_2+\mu}+\frac{x^2_3}{c^2_3+\mu}=1,\\
\frac{x^2_1}{c^2_1+\xi}+\frac{x^2_2}{c^2_2+\xi}+\frac{x^2_3}{c^2_3+\xi}=1,
\end{align*}
subject to the conditions $-c^2_3 < \xi < -c^2_2 < \mu < -c^2_1 < \rho$. Then the confocal ellipsoid $\p\GO$ is given by $\rho=\rho_0$ for some $\rho_0>0$.

Let
\beq\label{grho}
g(\rho) = (c_1^2+\rho)(c_2^2+\rho)(c_3^2+\rho),
\eeq
and define
\beq\label{vpform}
\Gvf_j (\rho)= \int_\rho^\infty \frac{1}{(c_j^2+s)\sqrt{g(s)}} ds, \quad j=1,2,3.
\eeq
Then the function $w$ defined by
\beq
w(\Bx) = \frac{1}{2} \int_\rho^\infty \frac{1}{\sqrt{g(s)}} ds - \frac{1}{2} \sum_{j=1}^3 \Gvf_j(\rho) x_j^2 + \frac{1}{2} \sum_{j=1}^3 \Gvf_j(\rho_0) x_j^2
\eeq
is a solution of \eqnref{free}. In fact, we can see that
$$
\pd{}{x_i} \left[\frac{1}{2} \int_\rho^\infty \frac{1}{\sqrt{g(s)}} ds - \frac{1}{2} \sum_{j=1}^3 \Gvf_j(\rho) x_j^2 \right] = \left( - \frac{1}{\sqrt{g(\rho)}} - \sum_{j=1}^3 \Gvf_j'(\rho) x_j^2 \right) \pd{\rho}{x_i} - \Gvf_i(\rho) x_i.
$$
Since
$$
\sum_{j=1}^3 \Gvf_j'(\rho) x_j^2 = - \sum_{j=1}^3 \frac{x_j^2}{(c_j^2+\rho)\sqrt{g(\rho)}} = -\frac{1}{\sqrt{g(\rho)}},
$$
we have
$$
\pd{}{x_i} \left[\frac{1}{2} \int_\rho^\infty \frac{1}{\sqrt{g(s)}} ds - \frac{1}{2} \sum_{j=1}^3 \Gvf_j(\rho) x_j^2 \right] =  - \Gvf_i(\rho) x_i,
$$
from which we see that
\beq
\nabla w(\Bx)= - (\Gvf_1(\rho)x_1, \Gvf_2(\rho)x_2, \Gvf_3(\rho)x_3) + (\Gvf_1(\rho_0)x_1, \Gvf_2(\rho_0)x_2, \Gvf_3(\rho_0)x_3).
\eeq
Using the relation
\beq\label{derivative}
\pd{\rho}{x_i}=\frac{2x_i}{c_i^2+\rho}\left[ \frac{x_1^2}{(c_1^2+\rho)^2}+\frac{x_2^2}{(c_2^2+\rho)^2}+\frac{x_3^2}{(c_3^2+\rho)^2}\right]^{-1},
\eeq
we see that $\GD w$ is constant. Note that $\nabla w=0$ on $\p\GO$ ($\rho=\rho_0$) and $\nabla w=\BA \Bx$ on $\p D$ where
\beq
\BA= \mbox{diag} [ \Gvf_1(\rho_0) - \Gvf_1(0), \Gvf_2(\rho_0) - \Gvf_2(0), \Gvf_3(\rho_0)-\Gvf_3(0) ].
\eeq
We emphasize that $\BA$ is negative-definite.

%%%%%%%%%%%%%%%%%%%%%%%%%%%%%%%%%%%%%%%%%%%%%%%%%%%%%%%
\section{Proof of Theorem \ref{conBall}}\label{sec:proof}
%%%%%%%%%%%%%%%%%%%%%%%%%%%%%%%%%%%%%%%%%%%%%%%%%%%%%%%

Let $w$ be the solution to \eqnref{free} with $\BA=c\BI$. We notice that $c \not=0$. Indeed, if $c = 0$, then we have
$$
0 \not=k|\Omega\setminus\overline{D}|= \int_{\Omega\setminus\overline{D}}\Delta w\ dx = \int_{\partial\Omega} \frac {\partial w}{\partial\nu}\ d\sigma - \int_{\partial D} \frac {\partial w}{\partial\nu}\ d\sigma = 0 -\int_{\partial D} \nu\cdot \Bd\ d\sigma = 0,
$$
which is a contradiction. Since $c\not=0$, by introducing new variables
$$
\By = \Bx + \frac 1c \Bd,
$$
we may assume that $\Bd =\bold 0$.  Set
\beq
A_{ij} = x_j \pd{}{x_i} - x_i \pd{}{x_j}, \quad i \neq j.
\eeq
It is worth mentioning that $A_{ij}$ is the angular derivative. Observe that $A_{ij}$ commutes with $\GD$, namely, $A_{ij} \GD = \GD A_{ij}$. So, we have $\GD A_{ij} w=0$ in $\GO \setminus \overline{D}$. Note that $A_{ij} w=0$ on $\p\GO$. Since $\nabla w(\Bx)= c\Bx$ on $\p D$, we see that $A_{ij} w=0$ on $\p D$. Then the maximum principle yields that
\beq
A_{ij} w=0 \quad\mbox{in } \GO \setminus \overline{D}.
\eeq
%So, $w$ is a radial function, in other words, $w(\Bx)=w(r)$ for $r=|\Bx|$.
Since $\GD w=k$ in $\GO \setminus \overline{D}$, $w$ satisfies the ordinary differential equation
\beq\label{ode1}
\frac{\p^2 w}{\p r^2} + \frac{2}{r} \frac{\p w}{\p r} =k \quad\mbox{in } \GO \setminus \overline{D}
\eeq
for $r=|\Bx|$. Choose a ball $B$ with $\overline{B} \subset \GO \setminus \overline{D}$. By \eqnref{ode1}, $w$ is of the form
\beq\label{ode2}
w(r)= \frac{k}{6} r^2 + \frac{k_1}{r} + k_2 \quad\mbox{in } \overline{B}
\eeq
for some real constants $k_1$ and $k_2$. Since $\Omega \setminus \overline{D}$ is connected and
$$
\GD\left(w- \frac{k}{6} r^2 - \frac{k_1}{r} - k_2\right)=0 \quad\mbox{in } \GO \setminus \overline{D},
$$
we have from \eqref{ode2}
\beq\label{ode3}
w(r)= \frac{k}{6} r^2 + \frac{k_1}{r} + k_2 \quad\mbox{in } \GO \setminus \overline{D}.
\eeq
Since $\pd{w}{r}=0$ on $\p\GO$, we must have
$$
  \frac{k}{3} r - \frac{k_1}{r^{2}}= 0 \quad\mbox{on } \partial\GO,
 $$
 and hence
 $$
 r^{3} =  \frac{3k_{1}}k\quad\mbox{on } \partial\GO.
 $$
 This means that $\partial\Omega = \partial B_{R}(\bold 0)$ for some $R >0$. Therefore we have

%Suppose that $\p\GO$ is not a sphere centered at the origin. Let $R$ be the smallest radius such that $\GO \subset B_R(0)$. Here $B_r(0)$ denotes the ball of radius $r$ centered at $0$. Then there is $\Ge>0$ such that for all $r$ satisfying $R-\Ge \le r \le R$, $\p B_r(0)$ intersect with $\p\GO$. Since $w=0$ on $\p\GO$, it means that $w(r)=0$ for $r \in [R-\Ge, R]$, which is not possible because of \eqnref{ode3}. So, $\p\GO= \p B_R(0)$.
$$
\nabla w(\Bx)= \frac{k}{3} \Bx - \frac{kR^3}{3} \frac{\Bx}{r^3}, \quad \Bx \in \GO \setminus \overline{D}.
$$
Since $\nabla w(\Bx)=c\Bx$ for all $\Bx \in \p D$, we must have
$$
\frac{k}{3} - \frac{kR^3}{3} \frac{1}{r^3} = c \quad\mbox{on } \partial D,
$$
or $r=\mbox{constant}$ for all $\Bx \in \p D$. It means that $\p D$ is a sphere centered at $\bold 0$. This completes the proof. \qed

%%%%%%%%%%%%%%%%%%%%%%%%%%%%%%%%%%%%%%%%%%%%%%
\section{Newtonian potential formulation}\label{sec:newton}
%%%%%%%%%%%%%%%%%%%%%%%%%%%%%%%%%%%%%%%%%%%%%%

In this section we reformulate the problem \eqnref{free} in terms of the Newtonian potentials and relate the problem with known characterization of ellipsoids using the property of the Newtonian potential.

Suppose that \eqnref{free} admits a solution and let $w$ be the solution.  Notice that by the second equation of \eqref{free}   $w$ is constant on each connected component of $\partial\Omega$, and by the third equation of \eqref{free}  $w(\Bx)= \frac{1}{2} \Bx \cdot \BA \Bx + \Bd\cdot\Bx + C$ for $\Bx \in \p D$ for some constant $C$.
Fix $\Bx \notin \overline{\GO} \setminus D$. We obtain from the divergence theorem that
\begin{align*}
k \int_{\GO \setminus D} \GG(\Bx-\By) d\By &= \int_{\GO \setminus D} \left[ \GD w(\By) \GG(\Bx-\By) - w(\By) \GD_\By \GG(\Bx-\By) \right] d\By \\
= -\int_{\p D} &\left[ \pd{w}{\nu}(\By) \GG(\Bx-\By) - w(\By) \pd{}{\nu_\By} \GG(\Bx-\By) \right] d\Gs(\By) - \int_{\p \Omega}  w(\By) \pd{}{\nu_\By} \GG(\Bx-\By)  d\Gs(\By)\\
= -\int_{\p D}& \left[ (\nu\cdot \BA \By +\nu\cdot\Bd) \GG(\Bx-\By) - (\frac{1}{2} \By \cdot \BA \By + \Bd\cdot\By +C) \pd{}{\nu_\By} \GG(\Bx-\By) \right] d\Gs(\By)\\
 &- \int_{\p \Omega}  w(\By) \pd{}{\nu_\By} \GG(\Bx-\By) d\Gs(\By).
\end{align*}

If $\Bx \in \Rbb^3 \setminus \overline{\GO}$, then  by dealing with the last integral  for each component of $\partial\Omega$ we have
$$
\int_{\p \Omega}  w(\By) \pd{}{\nu_\By} \GG(\Bx-\By) d\Gs(\By) = 0,
$$
and hence
\begin{align*}
k \int_{\GO \setminus D} \GG(\Bx-\By) d\By = - \mbox{Tr}\, \BA \int_{D} \GG(\Bx-\By) d\By.
\end{align*}
We can find a relation between $k$ and $\mbox{Tr}\, \BA$ from this formula. In fact, we have
$$
\lim_{|\Bx| \to \infty} \frac{1}{\GG(\Bx)} \int_{\GO \setminus D} \GG(\Bx-\By) d\By = |\GO \setminus D|
$$
and a similar formula for $\int_{D} \GG(\Bx-\By) d\By$. So we obtain
\beq
k |\GO \setminus D|= -\mbox{Tr}\, \BA |D|.
\eeq
We then see that
\beq
k \int_{\GO \setminus D} \GG(\Bx-\By) d\By + \mbox{Tr}\, \BA \int_{D} \GG(\Bx-\By) d\By = k|\GO| \left[ \widehat N_\GO(\Bx)- \widehat N_D(\Bx) \right],
\eeq
where $\widehat N_\GO$ and $\widehat N_D $ are the (averaged) Newtonian potentials on $\GO$ and $D$, respectively,  namely,
\beq
\widehat N_\GO(\Bx):= \frac{1}{|\GO|} \int_{\GO} \GG(\Bx-\By) d\By,
\eeq
and similarly for $\widehat N_D$.

If $\Bx \in D$, then we have for some constant $C^*$
\begin{align*}
k \int_{\GO \setminus D} \GG(\Bx-\By) d\By = - \mbox{Tr}\, \BA \int_{D} \GG(\Bx-\By) d\By + \frac{1}{2} \Bx \cdot \BA \Bx + \Bd \cdot \Bx + C^*,
\end{align*}
and hence
\beq
k|\GO| \left[ \widehat N_\GO(\Bx)- \widehat N_D(\Bx) \right] = \frac{1}{2} \Bx \cdot \BA \Bx + \Bd \cdot \Bx + C^*.
\eeq

We have shown that if \eqnref{free} admits a solution, then
\beq\label{dive}
\widehat N_\GO(\Bx)- \widehat N_D(\Bx) =
\begin{cases}
0, \quad & \Bx \in \Rbb^3 \setminus \GO, \\
\mbox{a quadratic polynomial}, \quad & \Bx \in D.
\end{cases}
\eeq

So we may reformulate the question: If \eqnref{dive} holds, then $D$ and $\GO$ are confocal ellipsoids. This is reminiscent of a question related to the Newton potential problem: If a Newtonian potential of a simply connected domain is a quadratic polynomial in the domain, then the domain must be an ellipsoid. This problem has been solved by Dive \cite{dive} and Nikliborc \cite{nikl} (see also \cite{DF} and \cite{km08}).

%%%%%%%%%%%%%%%%%%%%%%%%%%%%%%%%%%


\begin{thebibliography}{99}

\bibitem{book2} H. Ammari and H. Kang, \textsl{Polarization and moment
tensors with applications to inverse problems and effective medium
theory}, Applied Mathematical Sciences, Vol. 162, Springer-Verlag,
New York, 2007.

\bibitem{AKLL1} H. Ammari, H. Kang, H. Lee, and M. Lim, Enhancement of near cloaking using generalized polarization tensors vanishing structures. Part I: The conductivity problem, Comm. Math. Phys. 317 (2013), 253--266.

\bibitem{DF} E. DiBenedetto and A. Friedman, Bubble growth in porous media, Indiana Univ. Math. J. 35 (2) (1986), 573--606.

\bibitem{dive} P. Dive, Attraction des ellipsoides homog\`enes et r\'eciproques d'un
th\'eor\`eme de Newton, Bull. Soc. Math. France 59 (1931), 128--140.

\bibitem{glu} A. Greenleaf, M. Lassas, and G. Uhlmann, On nonuniqueness for Calderon's inverse problem,
Math. Res. Lett. 10 (2003), 685--693.

\bibitem{hashine} Z. Hashin, The elastic moduli of heterogeneous materials, J. Appl. Mech. 29 (1962), 143--150.

\bibitem{HS} Z. Hashin and S. Shtrikman, A variational approach to the theory of the effective magnetic permeability of multiphase materials, J. Appl. Phy. 33 (1962), 3125--3131.

\bibitem{JM} P. Jarczyk and V. Mityushev, Neutral coated inclusions of finite conductivity, Proc. R. Soc. A 468 (2012), 954--970.

\bibitem{kang} H. Kang, Conjectures of Polya-Szego and Eshelby, and the Newtonian potential problem; A review, Mechanics of Materials 41 (2009), 405--410.

\bibitem{kl14} H. Kang and H. Lee, Coated inclusions of finite conductivity neutral to multiple fields in two dimensional conductivity or anti-plane elasticity, Euro. J. Appl. Math., 25 (3) (2014), 329--338.

\bibitem{km08} H. Kang and G.W. Milton, Solutions to the P\'olya-Szeg\"o conjecture and the weak Eshelby conjecture, Arch. Ration. Mech. Anal. 188 (2008), 93--116.

\bibitem{kerker} M. Kerker, Invisible bodies, J. Opt. Soc. Am. 65 (1975), 376--379.

\bibitem{kohn1} R. V. Kohn, H. Shen, M. S. Vogelius, and M. I. Weinstein,
Cloaking via change of variables in electric impedance tomography,
Inverse Problems 24 (2008), article 015016.

\bibitem{milton}  G.W. Milton,
\newblock {\sl  The Theory of Composites},
\newblock Cambridge Monographs on Applied and Computational
Mathematics, Cambridge University Press, 2002.

\bibitem{MS} G. W. Milton and S. K. Serkov, Neutral coated inclusions in conductivity
and anti-plane elasticity, Proc. R. Soc. Lond. A 457 (2001), 1973--1997.

\bibitem{nikl} W. Nikliborc, Eine Bemerkung \"uber die Volumpotentiale, Math.
Zeit. 35 (1932), 625--631.

\bibitem{pendry} J. B. Pendry, D. Schurig, and D. R. Smith, Controlling electromagnetic fields, Science 312
(2006), 1780--1782.

\end{thebibliography}
\end{document}